\documentclass[opre,letterpaper,leqno,final]{informs3hack}

\SingleSpacedXI

\usepackage{placeins}
\usepackage[hypertexnames=false,colorlinks=true,breaklinks=true,bookmarks=true,urlcolor=blue,citecolor=blue,linkcolor=blue,bookmarksopen=false,draft=false]{hyperref}

\newtheorem{thm}{Theorem}[section]
\newtheorem{cor}[thm]{Corollary}
\newtheorem{lem}[thm]{Lemma}
\newtheorem{prop}[thm]{Proposition}
\newtheorem{defn}[thm]{Definition}
\newtheorem{rem}[thm]{Remark}

\DeclareMathOperator{\ldet}{ldet}
\DeclareMathOperator{\rank}{rank}
\DeclareMathOperator{\linx}{linx}

\DeclareMathOperator{\Diag}{Diag}
\DeclareMathOperator{\diag}{diag}

\TheoremsNumberedThrough     
\ECRepeatTheorems

\EquationsNumberedThrough    

\begin{document}

\RUNTITLE{Masking Anstreicher's linx bound}

\TITLE{\Large 
Masking Anstreicher's linx bound\\ for improved entropy bounds\thanks{M. Fampa was supported in part by CNPq grants 303898/2016-0 and
434683/2018-3. J. Lee was supported in part by AFOSR grant FA9550-19-1-0175.}}

\ARTICLEAUTHORS{%
\AUTHOR{Zhongzhu Chen}
\AFF{University of Michigan, Ann Arbor, \EMAIL{zhongzhc@umich.edu}} 
\AUTHOR{Marcia Fampa}
\AFF{Universidade Federal do Rio de Janeiro, Brazil, \EMAIL{fampa@cos.ufrj.br}}
\AUTHOR{Jon Lee}
\AFF{University of Michigan, Ann Arbor, \EMAIL{jonxlee@umich.edu}} 
} 

\RUNAUTHOR{Chen, Fampa, Lee}

\ABSTRACT{%
The maximum-entropy sampling problem is the NP-hard problem of maximizing the (log) determinant of 
an order-$s$ principle submatrix of a given order $n$ covariance matrix $C$. 
Exact algorithms are based on 
a branch-and-bound framework. The problem has wide applicability in spatial statistics, and in particular in 
environmental monitoring; see 
[C.-W. Ko, J. Lee, M. Queyranne  \emph{Oper. Res.,}
An  exact  algorithm  for  maxi-mum entropy sampling, 
43(4):684–691, 1995] and 
[J. Lee, Constrained   maximum-entropy sampling, \emph{Oper. Res.,}
 46(5):655–-664, 1998].
Probably the best upper bound for the maximum, empirically, is Anstreicher's scaled ``linx'' bound
(see [K.M.  Anstreicher.  Efficient  solution  of maximum-entropy sampling problems. \emph{Oper. Res.,} 68(6):1826--1835, 2020]). 
An earlier methodology for potentially improving any upper-bounding method is by masking; i.e. applying the
bounding method to $C\circ M$, where $M$ is any correlation matrix. 
We establish that 
the linx bound can be improved via masking by an amount that is at least  linear in $n$,
even when optimal scaling parameters are employed.
}%

\KEYWORDS{differential entropy, maximum-entropy sampling,
environmental monitoring, nonlinear combinatorial optimization, 
nonlinear integer optimization, convex relaxations}


\twocolumn 

\maketitle






\section{Introduction}
The \emph{maximum-entropy sampling problem}, a fundamental problem in optimal statistical design,
was formally introduced in \cite{SW} and then applied in many areas such as 
the re-design of environmental-monitoring networks (see \cite{ZidekSunLe2000}, for example). 
In the Gaussian case, the problem can be cast as
\begin{align*}
z(C,s)\!:=\! \max \{\ldet C[S,S] ~\!\!\!:~\!\!\!~ S\subseteq N, |S|=s\},\tag{MESP}\label{MESP}
\end{align*}
where ldet denotes the natural logarithm of the determinant, $C$ is an $n\times n$ covariance matrix (of Gaussian random variables), positive integer $s<n$ with $s\leq \rank(C)$, 
$N:=\{1,...,n\}$, and for subsets $S\subseteq N$ and $T\subseteq N$ ($1\le |S|\le n, 1\le |T|\le n$), $C[S,T]$ denotes the submatrix of $C$ having rows indexed by $S$ and columns indexed by $T$. 
Of course we should assume that $s\leq \rank(C)$, otherwise $z(C,s)=-\infty$.
We further assume that $C[j,j]>0$ for all $j\in N$,
because if we had any $C[j,j]=0$, then  such a $j$ 
could not be in any feasible solution of \ref{MESP}
having objective value greater than $-\infty$.
In the constrained version CMESP, we also have $m\geq 0$ side constraints: $\sum_{j\in S}a_{ij}\leq b_i$, for $1\leq i \leq m$. 

In the environmental-monitoring  application of \ref{MESP}, we collect time-series observations from  $n$ environmental monitoring stations, and we prepare a sample covariance matrix $C$ (see \cite{Al-ThaniLee1,Rpackage}, for details on how this can be done effectively).
In many situations, keeping all $n$ of the monitoring stations running is too costly, as so we wish to 
select a subset of size $s$, and continue monitoring only at them. Maximizing the 
``differential entropy'' of the $s$ sites, is a means of choosing the $s$-subset with maximum information.
In the U.S.A., for example, there are approximately $n=250$ stations that comprise the 
NADP (National Acidic Deposition Program) maintains the NTN (National Trends Network); see \cite{NADPNTN}).

 Unfortunately, \ref{MESP} is NP-hard (see \cite{KLQ}). Successful exact approaches are based on a branch-and-bound framework (see \cite{KLQ}). The original approach of \cite{KLQ} uses the so-called ``eigenvalue bound'' (the sum of the logs of the $s$ greatest eigenvalues of $C$)
 as an upper bound, and there has subsequently been considerable work on extensions and alternative bounding methods.
 \cite{LeeConstrained} extended the eigenvalue bound to CMESP, \cite{HLW} gave a ``combinatorial mask'' technique for improving the eigenvalue bound for \ref{MESP}, \cite{LeeWilliamsILP} further leveraged these
 ideas using integer-linear optimization, dynamic programming, and optimal matchings, and then
 \cite{AnstreicherLee_Masked,BurerLee} developed  the general ``masking'' technique to further improve 
 the eigenvalue bound. In another vein, bounding techniques based on convex-relaxation 
 for \ref{MESP} (which easily apply as
 well to CMESP), were developed; see \cite{AFLW_IPCO,AFLW_Using,Anstreicher_BQP_entropy,Kurt_linx}.
Although no bounding technique  wins on all instances, probably the current best is  the so-called ``linx bound''  of \cite{Kurt_linx}, which has several nice properties.
In recent work, \cite{chen_mixing} gave a general technique for 
 ``mixing'' more than one  bound that is based on convex relaxation. 
 
 \textbf{Notation.} 
$\Diag(x)\in \mathbb{R}^{n\times n}$ denotes the square diagonal matrix with diagonal elements as the elements of $x\in\mathbb{R}^n$. $\diag(X)\in\mathbb{R}^n$ denotes the vector obtained from the diagonal
elements of matrix $X\in  \mathbb{R}^{n\times n}$.
$\circ $ is the Hadamard (i.e., element-wise) product. 
$J_n$ is the square all-ones matrix of order $n$. 
$I_n$ is the identity matrix of order $n$. 
$\mathbf{e_n}$ is the $n$-vector of all-ones.
$\mathbf{e}^n_i$ denotes the $i$-th standard unit vector in $\mathbb{R}^n$. Often we will omit the size $n$ if it is clear from the context. 
If $S$ is a set, $|S|$ denotes the cardinality of the set. 


{\bf Scaling}. An important general technique for potentially improving some of the 
entropy upper bounds is ``scaling'',
based on the simple observation that 
for a positive constant $\gamma$, and $S$ 
with $|S|=s$, we have that 
\[
\det (\gamma C)[S,S] = \gamma^s \det C[S,S].
\]
With this identity,  we can easily see that
\[
z(C,s) = z(\gamma C,s) -s \log \gamma~.
\]
So upper bounds for $z(\gamma C,s)$
 yield upper bounds for $z(C,s)$,
shifting by $-s \log \gamma$. This idea was first exploited in \cite{AFLW_IPCO,AFLW_Using}.
It is worth noting that most bounding methods are \emph{not} invariant under
scaling; that is the bound does \emph{not} generally  shift by $-s\log \gamma$ (a notable exception being the eigenvalue bound).

{\bf Masking}. Another  important technique for potentially improving some of the 
entropy upper bounds is ``masking''.
Given a positive integer $n$, a \emph{mask} (also known as a ``correlation matrix'')
is an $n\times n$ symmetric positive-semidefinite  matrix
with $\diag(M)=\mathbf{e}$. We denote the set of order-$n$ masks as $\mathcal{M}_n$.
Masking for \ref{MESP}, introduced in full generality in \cite{AnstreicherLee_Masked}, 
is based on the  observation that for any $S\subset N$ and mask $M$,
\[
\det C[S,S] \leq \det (C\circ M)[S,S].
\]
That is, masking cannot decrease entropy. Therefore,
for any mask $M\in\mathcal{M}_n$, we have
\[
z(C,s) \leq z(C\circ M,s).
\]
\noindent 
This implies that upper bounds for $ z(C\circ M,s)$
are also upper bounds for $z(C,s)$. 
\medskip

{\bf The linx bound}.
\cite{Kurt_linx}  proposed the very-successful \emph{linx bound}
$\linx (C,s):= $
\begin{equation} \label{basiclinx}
\max\{ f(C,s;x) ~:~ x\in P(n,s)\},
\end{equation}
where
$f(C,s;x):= \frac{1}{2} \ldet(C\Diag(x)C+\Diag(\mathbf{e}-x))$, and 
$ P(n,s):=\{x\in \mathbb{R}^n ~:~ \mathbf{e}^\top x=s,~ 0\le x\le \mathbf{e}\}$.
Applying a scaling parameter $\gamma>0$, we arrive at the 
\emph{scaled linx bound} $\linx (C,s;\gamma):=$
\begin{equation} \label{gammalinx}
 \max\{ f(C,s;\gamma;x) ~\!:\!~ x\in P(n,s)\},
\end{equation}
where $f(C,s;\gamma;x):=$
\[
\textstyle \frac{1}{2} \left(\ldet(\gamma C  \Diag(x)C
+\Diag(\mathbf{e}-x))-s\log \gamma \right).
\]


We are interested in finding optimal scaling parameters in the context of the scaled linx bound.
Because \eqref{gammalinx} is an ``exact relaxation'' (i.e., the objective
of the relaxation on an $x\in\{0,1\}^n$ is exactly $\ldet C[S,S]$ for
$S$ equal to the support of $x$),
for every scaling parameter $\gamma>0$, 
the following useful fact (true for any exact relaxation)
is easy to see.
\begin{prop}\label{powersetopt}
If $\hat x\in \{0,1\}^n$ is an optimal solution of \eqref{gammalinx} for $\gamma=\hat \gamma$, then 
$\hat \gamma$ is optimal. That is,
$
	\linx  (C,s; \hat \gamma)=\min_{\gamma>0}\linx  (C,s; \gamma).
$
\end{prop}

An important result from  \cite{chen_mixing} is that by replacing the scaling parameter $\gamma$ with $e^{\psi}$, the linx bound  (as well as the  the BQP bound of  \cite{Anstreicher_BQP_entropy}) is convex in $\psi$, and so a locally-optimal scaling parameter is globally optimal (see \cite{chen_mixing} for  details and an algorithmic approach to optimizing the scaling parameter
based on this). 
In what follows, we further study the effect of the scaling parameter on the linx bound,
and we demonstrate the benefit masking can have on the linx bound.

Applying a mask $M\in\mathcal{M}_n$, we arrive at the \emph{scaled and masked linx bound} 
$\linx (C,s;M,\gamma):=$
\begin{equation} \label{maskgammalinx}
 \max\{ f(C,s;M,\gamma;x) ~\!:\!~ x\in P(n,s)\},
\end{equation}
where 
\begin{align*}
&f(C,s;M,\gamma;x)\! :=
\textstyle \!\frac{1}{2} \!
\left(\ldet(\!\gamma(C \!\circ\! M)\!\Diag(\!x\!)\!(C \!\circ\! M)\right.\\
&\left.\quad+\!\Diag(\mathbf{e}\!-\!x))\!-\!s\log \gamma \right).
\end{align*}
 We refer to the mask $M$ as  optimal for  $\linx (C,s;M,\gamma)$, if it minimizes its value among all  $M\in\mathcal{M}_n$.
For brevity, we write  $\linx (C,s; M) :=  \linx (C,s; M, 1)$ . 

We note that masking does not generally produce an exact relaxation, so Proposition \ref{powersetopt}
does not extend to a sufficient condition for optimal masks.

{\bf Main results and organization.} 
A main goal of ours is to demonstrate the strong potential for masking to improve the linx bound.
We do this by exhibiting  sequences $\{C_k,s_k;M_k,\gamma_k, \hat{\gamma}_k\}_{k=1}^\infty$, with $C_k\succeq 0$ of order $n_k$, $M_k\in\mathcal{M}_{n_k}$,
$n_k\to\infty$, such that $\linx(C_k,s_k;\gamma_k)-\linx(C_k,s_k;M_{k},\hat \gamma_k)\geq \alpha_k$, where $\alpha_k$ grows linearly with $n_k$. First we do this for $\gamma_k=\hat \gamma_k=1$ (i.e., no scaling). Then, at the expense of a worse lower bound $\alpha_k$,
we do this when $\gamma_k$ and $\hat \gamma_k$  are the optimal scale factors. 
To get such lower bounds on the gap $\linx(C_k,s_k;\gamma_k)-\linx(C_k,s_k;M_{n_k},\hat \gamma_k)$,
we need a good lower bound on $\linx(C_k,s_k;\gamma_k)$ and a good upper bound on $\linx(C_k,s_k;M_{n_k},\hat \gamma_k)$.
In fact, to establish these gaps, we will take $M_{n_k}=I_{n_k}$, and so to get our needed upper bounds,
we use an exact characterization of the linx bound and the optimal scaling for the linx bound,
when $C$ is diagonal (useful because $C_{n_k}\circ I_{n_k}$ is diagonal).
bn
In \S\ref{sec:gap}, we establish that using a mask but no scaling parameter (i.e., $\gamma=1$), the best-case improvement in the 
linx bound is at least linear in $n$; specifically, $\approx .0312n$.
In \S\ref{sec:gamma}, we study the  behavior of the linx bound
as we vary the scaling parameter $\gamma>0$. It was already established that
the linx bound is convex in $\log(\gamma)$ (see \cite{chen_mixing}). 
We establish the limiting behavior,
as $\gamma$ goes to 0 and to infinity. When $s=\rank(C)$, the limit  
as $\gamma$ goes to infinity can be better than any finite choice of $\gamma$;
in this case, we establish that the limit can be calculated by solving
a single convex optimization problem. In \S\ref{sec:gammagap}, we  establish that using a mask and \emph{optimal} scaling parameters, the best-case improvement in the 
linx bound remains at least linear in $n$; specifically, $\approx .024n$. \S\ref{sec:final} contains some final remarks.

\section{Linear gap for the linx bound}\label{sec:gap}

Our main goal  in this section is to establish a linear lower bound on the 
best-case gap between the linx bound and the masked linx bound, giving a good justification 
for considering mask optimization.
Specifically, we will give a sequence $\{C_n,s_n;M_n\}$, for all even positive integers $n$, with $C_n\succeq 0$ of order $n$, and $M_n\in\mathcal{M}_{n}$, such that $\linx(C_n,s_n)-\linx(C_n,s_n;M_{n})\geq \frac{1}{4}\log\left(\frac{4}{3}\right) n$. In fact, we will take $s_n:={\textstyle\frac{n}{2}}$, and $M_n:=I_n$. Because we use $M_n:=I_n$, we will 
have $\linx(C_n,s_n;M_{n})=\linx(\Diag(d_{(n)}),s_n)$, where $d_{(n)}=\diag(C_n)$. 
Hence it is useful to characterize, in general, the optimal solution of
\eqref{basiclinx} when $C$ is diagonal. Additionally, beyond our own use in the present work, we believe that such a 
characterization can be useful in future work on gaps for the linx bound.

Without loss of generality, we assume that $C:=\Diag(d)$ where $d\in\mathbb{R}^n$ and $d_1\!\geq\! d_2\!\geq\! \cdots\! \geq\! d_n\!>\!0$. Then
\[
f(C,s; x) = {\textstyle\frac{1}{2}} \log
\prod_{i\in N} 
\left(
(d_i^2-1)x_i+1
\right)
.
\]
\begin{lem} \label{lem:propxhat}
Let $C:=\Diag(d)$, where $d\in\mathbb{R}^n$ satisfies $d_1\geq d_2\geq \cdots \geq d_n > 0$. There exists   an optimal solution  $\hat x$ of \eqref{basiclinx} such that  $\hat x_1\ge \hat x_2 \ge \cdots \ge \hat x_n$ and $\hat x_i=\hat x_j$, for all $i,j\in N$, such that $d_i=d_j$.
\end{lem}

\proof{Proof.}
Clearly,  \eqref{basiclinx} has an optimal solution $\hat x$.  And
\begin{align*}
&\left((d_i^2-1)\hat x_i+1\right)\left((d_j^2-1)\hat x_j+1\right)\\
&-\!\left((d_i^2-1)\hat x_j+1\right)\!\left((d_j^2-1)\hat x_i+1\right)\\
&\!=\! (d_i^2-d_j^2)(\hat{x}_i-\hat{x}_j).
\end{align*}

If $d_i>d_j$, from the identity above we see that $\hat x_i\geq \hat x_j$, otherwise, by exchanging components $i$ and $j$ of $\hat x$, we would increase the objective value of \eqref{basiclinx}.

If $d_i=d_j$, let $\delta:= \hat x_i + \hat x_j $. Then,
\begin{align*}
&\left((d_i^2-1)\hat x_i+1\right)\left((d_j^2-1)\hat x_j+1\right)\\
& \quad = \left((d_i^2-1)\hat x_i+1\right)
 \left((d_j^2-1)(\delta -{ \hat x_i})+1\right)
\end{align*}
In this case, if $d_i=1$, the above function is constant. Otherwise, it is a univariate concave quadratic in $\hat{x_i}$, and its maximum is uniquely attained at $ \hat x_i =\delta/2$. Therefore, by setting $ \hat x_i = \hat x_j =\delta/2$, the maximum of the function is always attained. \Halmos
\endproof

\begin{defn} 
We refer to an optimal solution $\hat x$ of \eqref{basiclinx} which satisfies the properties in Lemma \ref{lem:propxhat} as a uniform  optimal solution.
\end{defn}

\begin{lem}\label{lem:conc}
Let $C:=\Diag(d)$, where $d\in\mathbb{R}^n$ satisfies $d_1\geq d_2\geq \cdots \geq d_n > 0$ and $0\leq x\leq \mathbf{e}$. Then $f(C,s;x)$ strictly increases with $x_i$, if $d_i>1$, does not change with $x_i$, if $d_i=1$, and strictly decreases with $x_i$, if $d_i<1$. Furthermore, $f(C,s;\cdot)$ is concave in $[0,1]^n$ and strictly concave if $d_i\neq 1$, for all $i\in N$.
\end{lem}

\proof{Proof.}
For all $i\in N$, $d_i>0$ and $0\leq x_i\leq 1$ implies that 
$(d_i^2-1)x_i+1  > 0$. 
Then, for $i\in N$,
\begin{align*}
&\frac{\partial f(C,s;x)}{\partial x_i}=\frac{d_i^2-1}{2 (
\!(d_i^2-1)x_i\!+\!1
)}\left\{\begin{array}{c}\!< 0, \mbox{ if }\! d_i<1,\\\!=0,\mbox{ if }\! d_i=1,\\\!>0,\mbox{ if }\! d_i>0,\end{array}\right.\\
&\frac{\partial^2 f(C,s;x)}{\partial x_i^2}=\frac{-(d_i^2-1)^2}{2 (\!
(d_i^2-1)x_i\!+\!1)^2}\left\{\begin{array}{c}\!< 0,\! \mbox{ if }\! d_i\neq 1,\\\!=0,\!\mbox{ if }\! d_i=1,\end{array}\right.\\
& \mbox{and}\\
&\frac{\partial^2 f(C,s;x)}{\partial x_i\partial x_j}=0, \mbox{ for } 1\leq i\neq j \leq n. \Halmos
\end{align*} 
\endproof

\begin{rem}
From Lemmas \ref{lem:propxhat} and \ref{lem:conc}, we see that if $C=\Diag(d)$, with $0<d_i\neq 1, \forall i\in N$, then  
\eqref{basiclinx} has a unique optimal solution, which is a uniform optimal solution.
\end{rem}


Next, we establish  necessary conditions  for $\hat x$ to be a uniform optimal solution for 
\eqref{basiclinx} when $C$ is diagonal,  based on
checking a  finite set of feasible directions; we could also
get these conditions from the KKT conditions for \eqref{basiclinx}, 
also establishing their sufficiency, but our approach is simpler and suits our purpose.

\begin{lem}\label{diagoptcondition}
Let $C:=\Diag(d)$, where $d\in\mathbb{R}^n$ satisfies $d_1\geq d_2\geq \cdots \geq d_n > 0$.
   Let $\hat x$ be a uniform  optimal solution of \eqref{basiclinx}. For $1\leq i<j\leq n$, we have
    \begin{equation}
    \frac{d_j^2-1}{(d_j^2-1)\hat x_j+1}\le \frac{d_i^2-1}{(d_i^2-1)\hat x_i+1}. \label{cond1}
\end{equation} 
Additionally, if $1 > \hat x_i \geq \hat x_j > 0$, then 
\begin{equation}
\frac{d_j^2-1}{(d_j^2-1)\hat x_j+1}= \frac{d_i^2-1}{(d_i^2-1)\hat x_i+1}. \label{cond2}
\end{equation}
\end{lem}

\proof{Proof.}
    \eqref{cond1} is clear when $\hat x_i=\hat x_j$, from the fact that $d_i\geq  d_j> 0$.
    So we may assume that $\hat x_i>\hat x_j$. In this case $\mathbf{e}_j-\mathbf{e}_i$ is a feasible direction 
    for $\hat{x}$ relative to  \eqref{basiclinx}. Because $\hat{x}$ is optimal for \eqref{basiclinx},
    we must have that $ \nabla f(C,s;\hat x)^\top (\mathbf{e}_j-\mathbf{e}_i)\leq 0$, which  is equivalent to \eqref{cond1}.

\eqref{cond2} follows from the fact that, in this case, $\mathbf{e}_i-\mathbf{e}_j$ is also a feasible 
direction  for $\hat{x}$ relative to \eqref{basiclinx}. \Halmos
\endproof
We have a corollary of Lemma \ref{diagoptcondition} for two special cases: $ d_n>1$ and  $d_1<1$. We will see later that the characterization of the optimal solution in general  can be reduced to the characterization of the optimal solution in these two special cases.
\begin{cor}\label{diagoptconditionspecial}
Let $C:=\Diag(d)$, where $d\in\mathbb{R}^n$ satisfies either $d_1\geq d_2\geq \cdots \geq d_n>1$ or $1>d_1\geq d_2\geq \cdots \geq d_n> 0$.
   Let  $\hat x$ be a uniform  optimal solution of \eqref{basiclinx}. Then,
    \begin{equation}\label{coreq1}
    \hat x_i-\hat x_j\le \frac{1}{d_j^2-1}-\frac{1}{d_i^2-1}.
\end{equation} 
Additionally, if $1 > \hat x_i \geq \hat x_j > 0$, then 
\begin{equation}\label{coreq2}
\hat x_i-\hat x_j=\frac{1}{d_j^2-1}-\frac{1}{d_i^2-1}.
\end{equation}
\end{cor}

\proof{Proof.}
    If either $d_n>1$ or $d_1<1$, we have $d_i^2-1\neq 0, \forall i\in N$. Also, both $d_i^2-1$ and  $d_j^2-1$ have the the same sign $\forall i, j\in N$. Together with $(d_i^2-1)\hat x_i+1>0, \forall i\in N$, we have that \eqref{cond1} and \eqref{cond2} equals \eqref{coreq1} and \eqref{coreq2}, respectively. \Halmos
\endproof

To characterize an optimal solution of \eqref{basiclinx} when $C$ is diagonal, we first establish a lemma that characterizes an optimal solution of \eqref{basiclinx} in the two special cases discussed in Corollary~\ref{diagoptconditionspecial}.

\begin{lem}\label{diagopt1}
Let $C:=\Diag(d)$, where $d\in\mathbb{R}^n$ satisfies either $d_1\geq d_2\geq \cdots \geq d_n >1$ or $1>d_1\geq d_2\geq \cdots \geq d_n > 0$.  Let  $\hat x$ be a uniform  optimal solution of \eqref{basiclinx}  for a given $0<s<n$. We have,
\begin{itemize}
    \item[(i)] if $ \frac{1}{d_{s+1}^2-1}-\frac{1}{d_s^2-1}\ge 1$,
     then  
    \[
    \hat x_i:=
  \left\{
  \begin{array}{ll}
1,&\hbox{for $1\leq i \leq s$,}\\
0,&\hbox{for $s+1\leq i \leq n$,}
  \end{array}
\right.
    \]
    \item[(ii)] if $\frac{1}{d_{s+1}^2-1}-\frac{1}{d_s^2-1} < 1$, 
     then $0<\hat x_s <1$, and
      \begin{equation}\label{eqlemma24}
    \hat x_i:=
  \left\{
  \begin{array}{ll}
\min\left\{1,~\hat x_s+\frac{1}{d_{s}^2-1}-\frac{1}{d_{i}^2-1}\right\},\\[6pt] 
 \qquad \qquad\qquad\hbox{for $1\leq i \leq s-1$,}\\
\max\left\{0,~\hat x_s+\frac{1}{d_{s}^2-1}-\frac{1}{d_{i}^2-1}\right\},\\[6pt] 
 \qquad \qquad\qquad \hbox{for $s+1\leq i \leq n.$}
  \end{array}
\right.
    \end{equation}
\end{itemize}
\end{lem}

\proof{Proof.}

We have already shown that under the hypotheses, \eqref{basiclinx} has a unique optimal solution $\hat x$, where $\hat x_1\ge \hat x_2\ge \cdots\ge \hat x_n$. Thus, for $(i)$, we only need to show that $\hat x_s=1$.
Suppose that $\hat x_s<1$; then 
$1>\hat x_s\ge \hat x_{s+1}>0$, i.e., $\hat x_{s}-\hat x_{s+1}<1\le \frac{1}{d_{s+1}^2-1}-\frac{1}{d_{s}^2-1}$, which violates the necessary condition \eqref{coreq2} in Corollary~\ref{diagoptconditionspecial}.

For $(ii)$, we see that if $\hat x_s=1$, then $\hat x_{s+1}=0$ and  $\textstyle \hat x_{s}-\hat x_{s+1}=1>\textstyle\frac{1}{d_{s+1}^2-1}-\textstyle \frac{1}{d_s^2-1}$, which violates the 
necessary 
condition \eqref{coreq1} in Corollary~\ref{diagoptconditionspecial}. If $\hat x_s=0$, then $\sum_{i=1}^n \hat x_i\le \sum_{i=1}^{s-1} \hat x_i\le s-1$, which contradicts the feasibility of $\hat{x}$. Therefore, $0<\hat x_s <1$. Finally, by the necessary 
conditions in Corollary~\ref{diagoptconditionspecial}, the other parts of $(ii)$ must hold. \Halmos
\endproof

In case $(ii)$ in Lemma~\ref{diagopt1}, 
we can solve  the equation $\mathbf{e}^\top x =s$  for $\hat x_s$:
\begin{align*}
&\sum_{i=1}^{s-1} \min\left\{1,~ \hat x_s+\frac{1}{d_{s}^2-1}-\frac{1}{d_{i}^2-1}\right\} +\hat x_s\\
&\quad 
+ \sum_{i=s+1}^n \max\left\{0,~ \hat x_s+\frac{1}{d_{s}^2-1}-\frac{1}{d_{i}^2-1}\right\}=s~,
\end{align*}
where $0<\hat x_s <1$. Note that the left-hand side of this equation is increasing, piecewise linear, and continuous in $\hat x_s$, so the equation is easy to solve. Once  $\hat x_s$ is determined, all $\hat x_i$, with $ i\neq s$ are also uniquely determined by \eqref{eqlemma24}.


Finally,  we have the following  characterization of optimal solutions when $C$ is diagonal.
\begin{thm}\label{diagoptall}
Let $C:=\Diag(d)$, where $d\in\mathbb{R}^n$ satisfies $d_1\!\geq\! d_2\!\geq\! \cdots\! \geq\! d_n \!> \!0$. Let $L:= \{i\in N : d_i <1\}$, $E:= \{i\in N : d_i =1\}$, and $G:= \{i\in N : d_i >1\}$. Then
$\hat x$ defined below is an optimal solution for \eqref{basiclinx}.  
\begin{itemize}
    \item[(i)] If $s\le |G|$, let $\tilde x\in \mathbb{R}^{|G|}$ be the optimal solution of 
$\widetilde{(1.1)}$,
which is \eqref{basiclinx}
   with $(C,s)$ replaced by $(\tilde{C},\tilde{s})$ as follows:
     $\tilde C:=\Diag(\tilde d)$, where $\tilde d\in \mathbb{R}^{|G|}, \tilde d_i:=d_i, 1\le i\le |G|$, and $\tilde s:=s$.
     Then,
    \[
    \hat x_i:=
  \left\{
  \begin{array}{ll}
\tilde x_i ,&\hbox{for $i\in G$,}\\
0,&\hbox{for $i\in E\cup L$.}
  \end{array}
\right.
    \]
\item[(ii)] If $|G|< s\le |G\cup E|$, let $\tilde x_i, i\in E$ be any value such that $0\le \tilde x_i\le 1$ and $\sum_{i\in E} \tilde x_i = s-|G|$. Then,
    \[
    \hat x_i:=
  \left\{
  \begin{array}{ll}
1 ,&\hbox{for $i\in G$,}\\
\tilde x_i,&\hbox{for $i\in E,$}\\
0 ,&\hbox{for $i\in L.$}
  \end{array}
\right.
    \]
\item[(iii)] If $|G\cup E|<s$, let $\tilde x\in \mathbb{R}^{|L|}$ be the optimal solution of 
 $\widetilde{(1.1)}$,
which is \eqref{basiclinx}
    with $(C,s)$ replaced by $(\tilde{C},\tilde{s})$ as follows:
$\tilde C:=\Diag(\tilde d)$ such that $\tilde d\in \mathbb{R}^{|L|}, \tilde d_i=d_{i+|G\cup E|}, 1\le i\le |L|$, and $\tilde s:=s-|G\cup E|$. Then,
    \[
    \hat x_i:=
  \left\{
  \begin{array}{ll}
1 ,&\hbox{for $i\in G\cup E$,}\\
\tilde x_i ,&\hbox{for $i\in L$.}
  \end{array}
\right.
    \]
\end{itemize}
\end{thm}

\proof{Proof.}

    We will prove $(i)$ in detail; $(ii)$ and $(iii)$ can be proved in a similar manner. The feasibility of $\hat x$ is obvious. We will prove that $\hat x$ is optimal as well. Let us assume otherwise, i.e, we assume that $x^*$ is an optimal solution  to \eqref{basiclinx} and \begin{equation}\label{hyp}f(C,s; x^*)> f(C,s;\hat x).\end{equation} 
    We first claim that  $x_i^*=0$, $\forall i\in E\cup L$. 
    Otherwise let  $x_i^*\! >\!0$, for some  $i \!\in\! E\cup L$. Then, as $s\le |G|$, by feasibility of $x^*$, we   have  $x_j^*\!< \!1$ for some $ j \in G$. Therefore $e_j -e_i $ is a feasible direction from $x^*$ in \eqref{basiclinx}. However, by Lemma \ref{lem:conc}, we have that   $\nabla f(C,s; x^*)^\top (e_j- e_i )>0$, contradicting the optimality of $x^*$. 
    
    Now, define
     $\tilde{x}^*\in \mathbb{R}^{|G|}$ such that $\tilde{x}^*_i=x^*_i$, $\forall i\in G$. As $x_i^*\!=\!0$, $\forall i\in E\cup L$, it is straightforward to see that  $\tilde{x}^*$ is feasible to $\widetilde{(1.1)}$. Thus $f(\tilde{C},\tilde{s};\tilde{x})\ge f(\tilde{C},\tilde{s};\tilde{x}^*)$. Note also that  $\hat x_i\!=\!0$, $\forall i\!\in\! E\cup L$, so $f(C,s;\hat{x})=f(\tilde{C},\tilde{s};\tilde{x})\ge f(\tilde{C},\tilde{s};\tilde{x}^*)=f(C,s;x^*)$, contradicting \eqref{hyp}.
\Halmos
\endproof

Having characterized an optimal solution for  \eqref{basiclinx} when $C $ is diagonal, we will now consider the more general case where $C$ is any positive-semidefinite matrix and establish a simple lower bound on $\linx (C,s)$, by considering  the eigenvalues of $C$.
\begin{lem} \label{lb}
 For any positive-semidefinite order-$n$ matrix $C$ and integer $0<s<n$, 
\begin{align*}
\linx (C,s)\ge {\textstyle\frac{1}{2}} \log
\prod_{i=1}^n \textstyle \left(\frac{s}{n}\lambda_i^2+1-\frac{s}{n}\right)
,
\end{align*}
where $\lambda_1,...,\lambda_n$ are the eigenvalues of $C$.
\end{lem}

\proof{Proof.}

We diagonalize $C$: That is, we choose an orthogonal matrix $Q$
so that $Q^\top C Q= \Lambda := \Diag (\lambda_1,\lambda_2,...,\lambda_n)$. 
	Let $\bar x=\frac{s}{n}\mathbf{e}$. Then 
	\begin{align*}
	\linx (C,s)\ge &f(C,s;\bar x)\\
	=&\textstyle\frac{1}{2}\ldet \left(\frac{s}{n}C^2+\left(1-\frac{s}{n}\right)I\right)\\
	=&\textstyle\frac{1}{2}\ldet \left(\frac{s}{n}Q\Lambda^2 Q^\top +\left(1-\frac{s}{n}\right)I\right)\\
	=&\textstyle\frac{1}{2}\ldet \left(\frac{s}{n}\Lambda^2+\left(1-\frac{s}{n}\right)I\right)\\
	=&{\textstyle\frac{1}{2}}\log
	\prod_{i=1}^n \textstyle \left(\frac{s}{n}\lambda_i^2+1-\frac{s}{n}\right)
	. \Halmos
	\end{align*}
\endproof

We now study the efficacy of using a mask $M$ for the linx bound (versus choosing $M=J$). We will show that there is a sequence of $\{C_n\}_{n\in \mathcal{I}}$ such that $\linx (C_n,{\textstyle\frac{n}{2}})-\linx (C_n,{\textstyle\frac{n}{2}}; I)\ge \frac{1}{4}\log(\frac{4}{3}) n\approx .0312 n$. The results shows that by choosing an appropriate mask $M$ different from $J$, we can decrease the linx bound by at least an amount that is linear in $n$. 

Recall that we have characterized an optimal solution of \eqref{basiclinx} when $C$ is diagonal  in Theorem~\ref{diagoptall} and a lower bound of $\linx(C,s)$ when $C$ is any positive-semidefinite matrix in Lemma~\ref{lb}. Note that  $C\circ I$ is diagonal. Then we have the following gap.
\begin{align}
&\linx (C,s)-\linx (C,s; I)\ge \label{maskgapineq}\\
&{\textstyle\frac{1}{2}}\!\log
\!\prod_{i=1}^n \!\!\left(\textstyle \frac{s}{n}\lambda_i^2\!+\!1\!-\!\frac{s}{n}\right)
\!-\! {\textstyle\frac{1}{2}}\!\log
\!\prod_{i=1}^n\!\!\left( d_i^2\hat x_i\!+\!1\!-\!\hat x_i\right)
\nonumber
\end{align}
where $\lambda_i, i\in N$ are the eigenvalues of $C$, $d_i, i\in N$ are diagonal elements of $C$, and $\hat x$ is an optimal solution of \eqref{basiclinx} with $C$ replaced by $C\circ I$. 
We will employ this lower bound on the gap in what follows.

Before presenting our main result, we will  characterize the optimal mask when $n=2$ and $s=1$.   
We will use this to construct a 
gap between $\linx (C,s)$ and $\linx (C,s; I)$ that is linear in the order of $C$. 

\begin{thm}\label{2orderoptmask}
    Let $C_2:=\left(\begin{array}{cc}
        a &    c  \\
        c&    b 
    \end{array}\right)$ be positive-semidefinite 
    where we assume, without loss of generality, $a\geq b$. Let 
    $M_2^*=\left(\begin{array}{cc}
         1 &m^*  \\
         m^*& 1 
    \end{array}\right)$ 
    be an optimal mask for  $\linx (C_2,1; M_2)$. We have,
\[
\begin{array}{cl}
    (i)&\!\! if\, c=0, \, then \, m^* \, can \, be\, any \,value\, in [-1,1];\\
    (ii)& \!\!if \;\frac{ab-1}{c^2}\ge 1, \;then \;m^*\; can \; be\;\pm 1;\\
    (iii)&\!\! if \;\frac{ab-1}{c^2}\le 0, \;then \;m^*\; can \; be\; 0;\\
    (iv)&\!\! if \;0<\frac{ab-1}{c^2}<1, \;then\; m^*\; can \; be\; \pm \sqrt{\frac{ab-1}{c^2}}.
\end{array}
\]
\end{thm}

\proof{Proof.}
    Let $M_2:=\left(\begin{array}{cc}
         1 &m  \\
         m& 1 
    \end{array}\right) \in \mathcal{M}_2$. Let $m^*=\argmin_{-1\le m\le 1}\left\{ \left(c^2m^2+1-ab\right)^2 \right\}$.
Considering that  $x_1+x_2=1$, we obtain
\begin{align*}
    &\ldet(\!(C_2\circ M_2)\!\Diag(\!x\!)\!(C_2\circ M_2)\!+\!I_2\!-\!\Diag(\!x\!)\!)\\
    =&\log (\!(c^2m^2+1-ab)^2 x_1 x_2 +(ax_1+bx_2)^2)\\
    \ge & \log (\!(c^2(m^*)^2+1-ab)^2 x_1 x_2 +(ax_1+bx_2)^2)\\
    =& \ldet(\!(C_2\circ M_2^*)\!\Diag(\!x\!)\!(C_2\circ M_2^*)\!+\!I_2\!-\!\Diag(\!x\!)\!),
\end{align*}
which implies that $M_2^*$ is an optimal mask. 

The values of $m^*$ in cases $(i-iv)$ can be easily obtained from $m^*=\argmin_{-1\le m\le 1}\left\{ \left(c^2m^2+1-ab\right)^2 \right\}$.
\Halmos
\endproof

For simplicity of the following discussions, we introduce the next lemma. 
\begin{lem}\label{2optgap}
 With the same hypotheses and notations as Theorem \ref{2orderoptmask}, $g(a,b,c):=\exp(2\linx(C_2,1;M_2^*))$. Define \[\Delta\linx(C,s;M):=\linx(C,s)-\linx (C,s; M).\] Then  
 \[
     \textstyle\Delta \linx (C_2,1; M_2^*)
     \ge\textstyle\frac{1}{2}\log\frac{(c^2+1-ab)^2+(a+b)^2}{4 g(a,b,c)}.
 \]
\end{lem}

\proof{Proof.}

Let $\lambda_1\ge \lambda_2\ge 0$ be the two eigenvalues of $C_2$. Considering that 
$\lambda_1+\lambda_2=a+b$ and $\lambda_1\lambda_2=ab-c^2$, the result follows directly from Lemma \ref{lb}.
\Halmos
\endproof

Note that for cases $(i-ii)$ in Theorem \ref{2orderoptmask}, there is no mask better than $J_2$. So we focus on  cases $(iii-iv)$. 
For  case $(iv)$, we can calculate from  the proof of Theorem \ref{2orderoptmask} that $\linx(C_2,1;M_2^*)=\frac{1}{2}\log\left(a^2\right)$. By Lemma \ref{2optgap}, we have 
\[
\textstyle\Delta\linx(C_2,1; M_2^*)
\ge\textstyle\frac{1}{2}\log\frac{(c^2+1-ab)^2+(a+b)^2}{4a^2}.
\]
 Note that $\frac{ab-1}{c^2}>0$ and $a\geq b$ imply $a>1$. 
Further, $a>1$, $\frac{ab-1}{c^2}<1$ and $ab\geq c^2$ 
 imply $0< c^2+1-ab\le 1 <a$. So,
\begin{align*}
    \textstyle\frac{1}{2}\log\frac{(c^2+1-ab)^2+(a+b)^2}{4a^2}< {\textstyle\frac{1}{2}}\log\left( {\textstyle\frac{5}{4}}\right).
\end{align*}
Moreover, by choosing $a=b=c>1$, we are in case $(iv)$, and the 
gap becomes $\frac{1}{2}\log(1+1/4a^2)$, which we can make as
close to $\frac{1}{2}\log(5/4)$ 
as we like.


For case $(iii)$, the optimal mask is $I_2$, and we can find a greater gap than we could for case $(iv)$. We prove this and our main result in the following theorem.

\begin{thm}\label{maskgap}
There is an infinite sequence of  positive-semidefinite matrices $\{C_n\}_{n\in \mathcal{I}}$ such that 
\begin{align*}
\linx \left(C_n,{\textstyle\frac{n}{2}}\right)\!-\!\linx \left(C_n,{\textstyle\frac{n}{2}}; I\right)\ge {\textstyle\frac{1}{4}}\log\left( {\textstyle\frac{4}{3}}\right) n,   \, n\in \mathcal{I},
\end{align*}
where $\mathcal{I}$ is the set of even integers. Moreover, for $n=2$, this is the maximum
possible lower bound on the gap that can be achieved using the lower bound from Lemma \ref{lb}.
%
\end{thm}

\begin{rem}
As we have indicated above in our analysis of case (iv), 
and proceeding similarly to 
how we proceed below, we can also get linear gaps with  masks that are
different from the identity mask, albeit with a worse constant (strictly less than $\frac{1}{4}\log(\frac{5}{4})$). 
\end{rem}

\proof{Proof.}(Theorem \ref{maskgap})
First, consider $n=2,s=1$. We use the same notations as in Theorem \ref{2orderoptmask} 
and  consider its case $(iii)$, where $ab\le 1$, so that the optimal mask is $I_2$. In the following, we will use $\hat x$ to denote the optimal solution of \eqref{maskgammalinx} for $C:=C_2$, $s=1$, $M:=I_2$, and $\gamma=1$;
so $\linx(C_2,1;I_2)=\linx(C_2,1;I_2,1)=f(C_2,1;I_2,1;\hat x)$.

We have two sub-cases to analyze:
\begin{itemize}
    \item[$(i)$] $a\ge 1\ge b$:  by Theorem~\ref{diagoptall},  $\hat x=(1,0)^\top$ is an optimal solution and $\linx(C_2,1;I_2)=$ $\frac{1}{2}\log\left(a^2\right)$. By Lemma \ref{2optgap}, we have
\begin{equation}\label{eqloga2}
\textstyle\Delta\linx(C_2,1;I_2) \ge \textstyle\frac{1}{2}\log\frac{(c^2+1-ab)^2+(a+b)^2}{4a^2}.
\end{equation}

Note that $(c^2+1-ab)^2+(a+b)^2\le 5a^2$, so $\frac{(c^2+1-ab)^2+(a+b)^2}{4a^2}\le \frac{5}{4}$. The equality can be obtained when $a=b=c=1$.
\item[$(ii)$] $1>a\ge b$: there are still two sub-cases: \[\textstyle \frac{1}{b^2-1}-\frac{1}{a^2-1}\ge 1 \mbox{ and } \frac{1}{b^2-1}-\frac{1}{a^2-1}< 1.\] 
\setlength{\leftmarginii}{0.0cm}
\begin{itemize}
    \item[$\bullet$] If $\frac{1}{b^2-1}-\frac{1}{a^2-1}\ge 1$, then by Theorem~\ref{diagoptall}, we also have $\hat x=(1,0)^\top$ and $\linx(C_2,1;I_2)=\frac{1}{2}\log\left(a^2\right)$. Thus, \eqref{eqloga2} also holds.
From $\frac{1}{b^2-1}-\frac{1}{a^2-1}\ge 1$ and $b^2\ge 0$, we can see that ${\textstyle\frac{1}{2}}\le a^2<1$ and $b\le \sqrt{2-\frac{1}{a^2}}$. Together with  $c^2\le ab<1$, letting $t:=\frac{1}{a^2}\in (1,2]$, we get
\begin{align*}
&\textstyle \max \frac{(c^2+1-ab)^2+(a+b)^2}{4a^2} =
\textstyle \max \frac{1+\left(a+\sqrt{2-\frac{1}{a^2}}\right)^2}{4a^2}\\
&\textstyle=\max \frac{1}{4} +\frac{3}{4} t-\frac{1}{4}t^2+{\textstyle\frac{1}{2}}\sqrt{2t-t^2}\\
& \textstyle\le \max \frac{1}{4} +\frac{3}{4} t-\frac{1}{4}t^2+\max {\textstyle\frac{1}{2}}\sqrt{2t-t^2}\\
& \textstyle=\frac{13}{16}+\frac{1}{2}< \frac{4}{3}.
\end{align*}
In the next sub-case, we will build a $\frac{1}{2}\log\left(\frac{4}{3}\right)$ gap, so the gap in the present sub-case is sub-optimal.
\item[$\bullet$] If $\frac{1}{b^2-1}-\frac{1}{a^2-1}<1$, then by Theorem~\ref{diagoptall}, 
$\hat x=\frac{1}{2}\left(1+\frac{1}{b^2-1}-\frac{1}{a^2-1},1-\frac{1}{b^2-1}+\frac{1}{a^2-1}\right)^\top$
is an optimal solution, and  we have
\begin{align*}
     \linx (C_2,1;& I_2)
    \!=\! \textstyle \frac{1}{2}\!\log\!\left(\!\frac{1}{4}\!\left(a^2\!+\frac{a^2-1}{b^2-1}\right)\!\left(b^2\!+\frac{b^2-1}{a^2-1}\right)\!\right).
\end{align*}
By Lemma \ref{2optgap}, we have
\[
\textstyle \Delta\linx (C_2,1; I_2) \ge \textstyle {\textstyle\frac{1}{2}}\log\frac{(c^2+1-ab)^2+(a+b)^2}{\left(a^2+\frac{a^2-1}{b^2-1}\right)\left(b^2+\frac{b^2-1}{a^2-1}\right)}.
\]
We claim that 
\[
\textstyle
\frac{(c^2+1-ab)^2+(a+b)^2}{\left(a^2+\frac{a^2-1}{b^2-1}\right)\left(b^2+\frac{b^2-1}{a^2-1}\right)}\leq \frac{1+(a+b)^2}{\left(a^2+\frac{a^2-1}{b^2-1}\right)\left(b^2+\frac{b^2-1}{a^2-1}\right)}\leq \frac{4}{3}.
\]
The first inequality holds because $c^2\leq ab<1$ and the second  holds for being equivalent to 
\[
\textstyle (1\!-\!2ab)^2+(a\!-\!b)^2+4(a^2\!-\!b^2)\left(\frac{1}{b^2-1}-\frac{1}{a^2-1}\right) \geq 0,
\]
We get equality in both with $a=b=c=\frac{\sqrt{2}}{2}$.
\end{itemize}
\end{itemize}

\vspace{0.1in}
\noindent
In the analysis above, we  see that we can create the largest gap in the last case. Therefore, we define \[\textstyle C_2:=\left(\begin{array}{cc}
       \frac{\sqrt{2}}{2}  & \frac{\sqrt{2}}{2} \\
       \frac{\sqrt{2}}{2}  & \frac{\sqrt{2}}{2}
    \end{array}\right).\]  Then $\lambda_1=\sqrt 2,\lambda_2=0$, the optimal solution for $C_2\circ I_2$ is $\left({\textstyle\frac{1}{2}},{\textstyle\frac{1}{2}}\right)^\top$, and 
    \begin{align*}
        &\linx (C_2,1)-\linx (C_2,1; I_2)\\
        &\ge {\textstyle\frac{1}{2}}\log\left(\left({\textstyle\frac{1}{2}}\cdot (\sqrt 2)^2+1-{\textstyle\frac{1}{2}}\right)\left({\textstyle\frac{1}{2}}\cdot 0^2+1-{\textstyle\frac{1}{2}}\right)\right)\\
        &\textstyle
        -{\textstyle\frac{1}{2}}\log\left(\!\!\left({\textstyle\frac{1}{2}}\cdot \left(\!\frac{\sqrt{2}}{2}\right)^2\!+\!1-{\textstyle\frac{1}{2}}\!\right)\!\!\!\left(\!{\textstyle\frac{1}{2}}\cdot \left(\frac{\sqrt{2}}{2}\right)^2\!+\!1-{\textstyle\frac{1}{2}}\!\right)\!\!\right)\\
       & =\textstyle {\textstyle\frac{1}{2}}\log\left(\frac{4}{3}\right).
    \end{align*}
For $n=2k$, we construct a block-diagonal matrix $C_n$ with $k={\textstyle\frac{n}{2}}$ blocks, and each block is such a $C_2$ matrix. Then we take $s={\textstyle\frac{n}{2}}$. In this way, $C_n$ has $k$ eigenvalues of $\sqrt{2}$ and $k$ eigenvalues of $0$. Also, all diagonal elements of $C_n$ are $\frac{\sqrt{2}}{2}$.  By \eqref{maskgapineq}, we have
\begin{align*}
&\linx \left(C_n,{\textstyle\frac{n}{2}}\right)-\linx \left(C_n,{\textstyle\frac{n}{2}}; I\right)\\
 & \textstyle  \ge  \frac{n}{4}\log\left(\left({\textstyle\frac{1}{2}} (\sqrt 2)^2+1-{\textstyle\frac{1}{2}}\right)\left({\textstyle\frac{1}{2}} (0)^2+1-{\textstyle\frac{1}{2}}\right)\right)\\
& ~ \textstyle -\frac{n}{4}\!\log\!\left(\!\!\left(\!{\textstyle\frac{1}{2}} \left(\frac{\sqrt{2}}{2}\right)^2+1-{\textstyle\frac{1}{2}}\right)\!\!\!\left(\!{\textstyle\frac{1}{2}} \left(\frac{\sqrt{2}}{2}\right)^2+1-{\textstyle\frac{1}{2}}\right)\!\!\right)\\
& 
=\textstyle \frac{1}{4}\log\left(\frac{4}{3}\right)n. \Halmos
\end{align*}

\endproof


\section{Optimal scaling parameter:  some special cases and general behavior}\label{sec:gamma}
In this section, we first show how an appropriate scaling parameter $\gamma$ can help improve the linx bound by forcing one optimal solution of \eqref{gammalinx} to lie in $\{0,1\}^n$ when $C$ is diagonal or $C$ is non-singular of order $2$. Next, we show the following results: $(i)$ if $s<\rank(C)$, then an optimal scaling parameter $\gamma$ for \eqref{gammalinx}  can always  be obtained, $(ii)$ if $s=\rank(C)$ and $\hat \gamma$ is an optimal scaling parameter for linx, then so is any $\gamma\geq \hat\gamma$, $(iii)$ if $s>\rank(C)$, there is no optimal $\gamma$. In fact, in this case we show that the linx-bound has the nice property of recognizing the behavior of \ref{MESP}, it tends to minus infinity as $\gamma$ tends to infinity.

\begin{prop}\label{diagoptgamma}
For diagonal positive-definite matrix $C:=\Diag \{d_1,...,d_n\}$, where $d_1\ge ...\ge d_n> 0$ and $0<s<n$, the scaling parameter $\hat \gamma=\frac{1}{d_s^2}$ forces an optimal solution of \eqref{gammalinx} to lie in $\{0,1\}^n$. 
	Therefore $\hat \gamma$ is an optimal scaling parameter.
\end{prop}

\proof{Proof.}

	Note that
	\[
	f(C,s; \gamma; x)={\textstyle\frac{1}{2}} \log
	\prod_{i=1}^n (\gamma d_i^2x_i+1-x_i)
	-{\textstyle\frac{1}{2}} s\log\gamma.
	\]
Partition $N$ as $N=L'\cup E'\cup G'$, where $\gamma d_i^2<1,i\in L'$;  $\gamma d_i^2=1, i\in E'$; $\gamma d_i^2>1, i\in G'$. As we have seen in Lemma \ref{lem:conc}, $f(C,s;  \gamma; x)$ strictly decreases with $x_i, i\in L'$, does not change with $x_i, i\in E'$ and strictly increases with $x_i, i\in G'$. So, if there is a $\gamma>0$ such that $|G'|\le s$ while $|E'\cup G'|\ge s$. By Theorem \ref{diagoptall},
	$(e_s^\top ,0)^\top $ is an optimal solution for \eqref{gammalinx} which lies in $\{0,1\}^n$. In fact, $\hat \gamma:=\frac{1}{d_s^2}$ is  such a  
	scaling 
	parameter. Therefore, by Proposition ~\ref{powersetopt}, $\hat \gamma$ is  optimal.
\Halmos
\endproof

\begin{prop}\label{2optgamma}
Let $C_2:=\left(\begin{array}{cc}
        a &    c  \\
        c&    b 
    \end{array}\right)$ be positive-definite 
    where we assume, without loss of generality, $a\geq b$. Let $s=1$.
     Then the scaling parameter $\hat \gamma=\frac{a^2-c^2}{(ab-c^2)^2}$ forces an optimal solution of \eqref{gammalinx} to lie in $\{0,1\}^2$. Therefore $\hat \gamma$ is an optimal scaling parameter.
\end{prop}

\proof{Proof.}

We have
	\begin{align*}
	&f(C_2,1; x)\\
	&={\textstyle\frac{1}{2}} \ldet (C_2 \Diag (x) C_2+I_2-\Diag (x)) \\ 
	&={\textstyle\frac{1}{2}} \!\log (\!(c^2+1-ab)^2 x_1 x_2 +(ax_1+bx_2)^2).
	\end{align*}
	Because $f(C_2,1;x)$ is concave, and the null space of $\mathbf{e}_2$ is $\{(t,-t)^\top ~:~ t\in \mathbb{R}\}$, to prove that one optimal solution lies in $\{0,1\}^2$ (in particular, we assume this optimal solution is $\hat x=(1,0)^\top$), we only need to prove 
	\[
	\left.\frac{f(C_2,1; \hat x - t(e_1-e_2))}{\partial t}\right|_{t=0} \le 0
	\]
	which is equivalent to
\begin{equation}\label{eq:diff}
\frac{\partial f(C_2,1; \hat x)}{\partial x_1}- \frac{\partial f(C_2,1; \hat x)}{\partial x_2}\ge 0,
\end{equation}
and finally
$
2(a^2-c^2)-(c^2-ab)^2-1\ge 0.
$

Because $f(C_2,s; \gamma; x)=f(\sqrt{\gamma} C_2,s; x)+{\textstyle\frac{1}{2}}\log\gamma$, we have that 
\[
\frac{\partial f(C_2,1;\gamma; \hat x)}{\partial x_1}- \frac{\partial f(C_2,1; \gamma;\hat x)}{\partial x_2}\ge 0,
\] 
is equivalent to
\[
\frac{\partial f(\sqrt{\gamma}C_2,1; \hat x)}{\partial x_1}- \frac{\partial f(\sqrt{\gamma}C_2,1; \hat x)}{\partial x_2}\ge 0,
\]
and finally 
\begin{equation} \label{2ordersuffcondi}
2(a^2-c^2)\gamma-(c^2-ab)^2\gamma^2-1\ge 0.
\end{equation}
    The left-hand side of \eqref{2ordersuffcondi} is maximized by $\hat\gamma=\frac{a^2-c^2}{(ab-c^2)^2}$ and the corresponding value is $\frac{(a^2-c^2)^2}{(ab-c^2)^2}-1$ which is nonnegative because $a^2\ge ab>c^2$. Note that if $a>b$ then $\frac{(a^2-c^2)^2}{(ab-c^2)^2}-1>0$, which means there is an $\epsilon>0$ such that for any $\gamma\in [\hat\gamma-\epsilon,\hat\gamma+\epsilon]$, $\frac{(a^2-c^2)^2}{(ab-c^2)^2}-1\ge 0$ and $\hat x$ is an optimal solution, i.e., the optimal scaling parameter is not unique.  Finally, by Proposition \ref{powersetopt}, $\hat \gamma$ is optimal. \Halmos
    \endproof
    

Not unexpectedly, there also exists a large and simple class of $C$ where no optimal solution of \eqref{gammalinx} lies in $\{0,1\}^n$ for any $0<s<n$ and any scaling parameter $\gamma$, as we will see in Theorem \ref{neverinpowerset}.

First, note that when $C=\tau_1 I$, for any $\tau_1>0$, in all cases of Theorem~\ref{diagoptall}, $\hat x:= \frac{s}{n} \mathbf{e}$ is an optimal solution of \eqref{basiclinx}. The same observation can be extracted  from Lemma \ref{lem:propxhat}. In fact, this observation is also a special case  of the following result, which follows immediately from the concavity of $f(C,s;x)$ and its symmetry in this case.
 
\begin{prop}
    \label{tauijopt}
	Suppose that $\tau_1>0$, $\tau_2 \ge 0$,and $0<s<n$ integer.
	Let $C=\tau_1 I+\tau _2 J$, then $\hat x=\frac{s}{n}\mathbf{e}$ is an optimal solution for \eqref{basiclinx}.
\end{prop}


\begin{thm}\label{neverinpowerset}
	For any order $n\ge 3$, any $0<s<n$  and $C=\tau_1 I+\tau_2 J$, $\tau_1>0$, $\tau_2> 0$, and for any scaling parameter $\gamma>0$, the optimal solution of \eqref{gammalinx} cannot lie in $\{0,1\}^n$.
\end{thm}

\proof{Proof.}

	By Proposition~\ref{tauijopt}, one optimal solution for \eqref{gammalinx} is $\frac{s}{n} \mathbf{e}$ under the setting of this theorem. From the proof of [Theorem 21, \cite{chen_mixing}], we have that if 
	\begin{align}
&	\gamma C\Diag (y)C-\Diag (y)\neq 0,\label{strictconcavecondi}\\
&	\hbox{	when $\mathbf{e}^\top y=0$, $-\mathbf{e}\le y\le \mathbf{e}$, $y\neq 0$}, \nonumber
	\end{align}
	then $f(C,s;\gamma;x)$ is strictly concave with a unique optimal solution on the feasible region of \eqref{gammalinx}. Because $\frac{s}{n}\mathbf{e}$ is already optimal in this case, we see that the optimal solution cannot lie in $\{0,1\}^n$. Now we prove that \eqref{strictconcavecondi} holds.
	
	Substituting $\tau_1 I+\tau_2 J$  for $C$ in \eqref{strictconcavecondi}  and dividing by $\gamma$, we get 
	\begin{align} 
	& \textstyle
	\left(\tau_1^2-\frac{1}{\gamma}\right)\Diag (y)+\tau_1\tau_2(\mathbf{e}y^\top +y\mathbf{e}^\top )\neq 0 \label{strictconcavecondi2}
	\end{align}
It is easy to see that if \eqref{strictconcavecondi2} is \emph{not} satisfied, then $y_i+y_j=0$ for all $i\not=j$.
But  this cannot be true when $n\geq 3$ for $y\not=0$. \Halmos
\endproof


Interestingly, there is also  a very simple example for which no $\gamma>0$ can be an optimal scale factor :
\begin{prop}\label{2nooptgamma}
	For $C:=J_2$, $s:=1$, there is no optimal scaling factor $\gamma$ for \eqref{gammalinx}. In fact, for all $\gamma>0$,
	\[
	\textstyle
	\linx (J_2,1; \gamma)={\textstyle\frac{1}{2}}\log \left(1+\frac{1}{4\gamma}\right).
	\]
	which monotonically decreases as $\gamma$ increases.
\end{prop}

\proof{Proof.}

	\begin{align*} 
& 	\linx (J_2,1;\gamma)\\
& =  {\textstyle\frac{1}{2}}\max _{\begin{array}{c}
	  \scriptstyle   x_1+x_2=1 \\
	  \scriptstyle   0\le x_1, x_2\le 1
	\end{array}}\left\{\log\left(\gamma (x_2+x_1)(2-x_1-x_2)\right.\right.\\[-25pt]
	&\qquad\qquad\qquad\qquad \left. \left.
	+(1-x_1)(1-x_2)\right)-\log \gamma
	\right\}\\[5pt]
&	\le  {\textstyle\frac{1}{2}}\max_{\begin{array}{c}
	    \scriptstyle   x_1+x_2=1 \\
	   \scriptstyle    0\le x_1, x_2\le 1
	\end{array}}\left\{ 
	\log\left(\vphantom{\left(\textstyle\frac{2-x_1-x_2}{2}\right)^2}\right.\gamma (x_2+x_1)(2-x_1-x_2)\right.\\[-20pt]
	&\qquad\qquad\qquad\qquad\qquad\left.\left. +\left(\textstyle\frac{2-x_1-x_2}{2}\right)^2\right)-\log \gamma
	\right\}\\
&\textstyle	=  \textstyle\frac{1}{2}\left(\log \left(\gamma+\frac{1}{4}\right)-\log\gamma\right)	= {\textstyle\frac{1}{2}}\log \left(1+\frac{1}{4\gamma}\right).
	\end{align*}
	
Note that both maximums are achieved at $x_1=x_2=1/2$, and so the inequality is an equation.
\Halmos
\endproof

Based on Proposition~\ref{2nooptgamma}, our interest is in what cases, we are guaranteed to have a finite optimal scaling parameter $\gamma$. In fact, a broad sufficient condition is $s<\rank (C)$ by the following theorem.

\begin{thm}\label{ssmallrank}
	For all positive-semidefinite $C$ and $0<s<n$, we have
	$$
	\lim_{\gamma\rightarrow 0} \linx  (C,s; \gamma)=+\infty.
	$$
	If we further assume that $s<\rank(C)$, then
	$$
	\lim_{\gamma\rightarrow +\infty} \linx  (C,s; \gamma)=+\infty.
	$$
\end{thm}

 \proof{Proof.}
 
	For all $\gamma>0$, by setting $\bar x:=\frac{s}{n} \mathbf{e}$, we have
	\begin{align*}
	&\linx (C,s; \gamma)\ge  f(C,s; \gamma; \bar x)\\
	& \textstyle \quad ={\textstyle\frac{1}{2}} \left(\ldet \left(\gamma\frac{s}{n}C^2+I-\frac{s}{n}I\right)-s\log \gamma\right).
	\end{align*}
	When $\gamma\rightarrow 0$, $\gamma\frac{s}{n}C^2+I-\frac{s}{n}I\rightarrow \left(1-\frac{s}{n}\right)I$. So
	\begin{align*}
	&\textstyle \lim_{\gamma\rightarrow 0} \ldet \left(\gamma\frac{s}{n}C^2+I-\frac{s}{n}I\right)= n\log \left(1-\frac{s}{n}\right), \\
	&\hbox{and } \lim_{\gamma\rightarrow 0}-s\log\gamma=+\infty.
	\end{align*}
	Therefore, $\lim_{\gamma\rightarrow 0}\linx (C,s; \gamma)= +\infty$.

But we can also write $ f(C,s; \gamma; \bar x)=$
\[
 \textstyle \frac{1}{2}\left(\ldet \left(\frac{s}{n}C^2+\frac{1}{\gamma}\left(1-\frac{s}{n}\right)I\right)\!+\!(n\!-\!s)\log\gamma\right).
 \]
Note that
$\lim_{\gamma\rightarrow +\infty}(n-s)\log\gamma=+\infty$.  Further, 
if $C$ is non-singular, then
	\begin{align*}
	\lim_{\gamma\rightarrow +\infty} 	\textstyle\ldet \left(\frac{s}{n}C^2+\frac{1}{\gamma}\left(1-\frac{s}{n}\right)I\right)= \ldet \left(\frac{s}{n}C^2\right).
	\end{align*}
So we can conclude that $\lim_{\gamma\rightarrow +\infty} \linx (C,s; \gamma)=+\infty$, when $C$ is nonsingular.

	

	If $C$ is singular, then we have
	\begin{align*}
		\textstyle
	\lim_{\gamma\rightarrow +\infty}\ldet \left(\frac{s}{n}C^2+\frac{1}{\gamma}\left(1-\frac{s}{n}\right)I\right)= -\infty,
	\end{align*}
	and we cannot immediately conclude anything useful.
	So we proceed differently.
	When $s<\rank (C)$, without loss of generality,
	we can write $C=Q\Lambda Q^\top$, where $Q$ is orthogonal and $\Lambda:=\Diag(\lambda_1,...,\lambda_n)$ with $\lambda_1\ge \lambda_2\ge ...\ge \lambda_n\ge 0$. We have $\lambda_i\neq 0$ for $i\le \rank (C)$ and $\lambda_i=0$ for $i> \rank (C)$. 
	By 
	L'H\^optital's
	rule,
	\begin{align*}
	&\	 \lim_{\gamma\rightarrow +\infty} \textstyle\frac{\ldet \left(\frac{s}{n}C^2+\frac{1}{\gamma}\left(1-\frac{s}{n}\right)I\right)}{(n-s)\log\gamma}\\
	=&\lim_{\gamma\rightarrow +\infty} \textstyle\frac{\partial \left(\ldet \left(\frac{s}{n}C^2+\frac{1}{\gamma}\left(1-\frac{s}{n}\right)I\right)\right)/\partial \gamma}{\partial \left((n-s)\log\gamma\right) /\partial \gamma}\\
	=& \lim_{\gamma\rightarrow +\infty}\textstyle\frac{\text{tr}\left(\left(\frac{s}{n}C^2+\frac{1}{\gamma}\left(1-\frac{s}{n}\right)I\right)^{-1}\left(1-\frac{s}{n}\right)I\right)\frac{-1}{\gamma^2}}{(n-s)\frac{1}{\gamma}} \\
	=& \lim_{\gamma\rightarrow +\infty}\textstyle\frac{-1}{n\gamma} \text{tr}\left(\left(\frac{s}{n}C^2+\frac{1}{\gamma}\left(1-\frac{s}{n}\right)I\right)^{-1}\right)\\
	=& \lim_{\gamma\rightarrow +\infty}\textstyle\frac{-1}{n} \text{tr}\left(\left(\gamma\frac{s}{n}C^2+\left(1-\frac{s}{n}\right)I\right)^{-1}\right)\\
	=& \lim_{\gamma\rightarrow +\infty}\textstyle\frac{-1}{n} \text{tr}\left(\left(\gamma\frac{s}{n}Q\Lambda^2Q^\top+\left(1-\frac{s}{n}\right)QQ^\top\right)^{-1}\right)\\
	=&\lim_{\gamma\rightarrow +\infty}\textstyle\frac{-1}{n} \text{tr}\left(Q\left(\gamma\frac{s}{n}\Lambda^2+\left(1-\frac{s}{n}\right)I\right)^{-1}Q^\top\right)\\
	=&\textstyle\frac{-1}{n-s}\text{tr}\left(\text{Diag}\left\{\mathbf{0}_{1\times \rank (C)},\mathbf{e}_{n-\rank (C)}^\top\right\}\right)\\
	=&\textstyle-\frac{n-\rank (C)}{n-s}.
	\end{align*}
	This means for every $\epsilon>0$, there exists $\gamma_{\epsilon}>0$ such that when $\gamma>\max\{\gamma_{\epsilon},1\}$,
	\begin{align*}
	&\textstyle
	 {\textstyle\frac{1}{2}}\left(\ldet \left(\frac{s}{n}C^2+\frac{1}{\gamma}\left(1-\frac{s}{n}\right)I\right)+(n-s)\log\gamma\right)\\
	& \textstyle\quad \ge  \frac{1}{2}\left(-\frac{n-\rank (C)}{n-s}-\epsilon+1\right)(n-s)\log\gamma.
	\end{align*}
	So 
	
	\vspace{-20pt}
	\begin{align*}
	&\lim_{\gamma\rightarrow +\infty} \linx (C,s; \gamma)\\
	\ge & \lim_{\epsilon\rightarrow 0}\lim_{\gamma\rightarrow +\infty}\textstyle\left(-\frac{n-\rank (C)}{n-s}-\epsilon+1\right)(n-s)\log\gamma\\
	= & \lim_{\epsilon\rightarrow 0}\lim_{\gamma\rightarrow +\infty}\textstyle\left(\frac{\rank (C)-s}{n-s}-\epsilon\right)(n-s)\log\gamma\\
	=& +\infty. \Halmos
	\end{align*}
\endproof

\begin{cor}
For all positive-semidefinite $C$ and $0<s<n$ where $s<\rank(C)$, we can find a finite optimal scaling parameter $\hat \gamma$ such that
\begin{align*}
\linx (C,s;\hat\gamma)=\min_{\gamma>0} \linx (C,s; \gamma).
\end{align*}
\end{cor}

\proof{Proof.}

By \cite{chen_mixing}, if we replace $\gamma$ with $e^{\psi}$, then $\linx (C,s; e^{\psi})$ is convex in $\psi$ and by Theorem~\ref{ssmallrank}, 
\begin{align*}
\lim_{\psi\rightarrow -\infty}  \linx (C,s; e^{\psi})=\lim_{\gamma\rightarrow 0} \linx  (C,s; \gamma)&=+\infty\\
\lim_{\psi\rightarrow +\infty}  \linx (C,s;  e^{\psi})=\lim_{\gamma\rightarrow +\infty} \linx  (C,s; \gamma)&=+\infty.
\end{align*}
We can conclude that a minimizing $\hat{\psi}$ exists, and then we have the 
minimizer $\hat \gamma:=e^{\hat{\psi}}$. 
\Halmos
\endproof

When $s=\rank(C)$, the following result establishes that  $\lim_{\gamma\rightarrow \infty}\linx (C,s; \gamma)$ exists and is finite, and $\linx (C,s; \gamma)$ is monotonically non-increasing in $\gamma$. This implies that if $\hat \gamma>0$ is optimal, then all $\gamma\ge \hat \gamma$ are optimal.

\begin{thm}\label{seqrank}
	When $s=\rank(C)$, without loss of generality, we can write $C$ as $C=Q\Lambda Q^\top $, where $Q$ is orthogonal and $\Lambda: =\Diag (\lambda_1,...,\lambda_n)$ with $\lambda_1\ge...\ge\lambda_s> \lambda_{s+1}=... =\lambda_n=0$. Denote $\Lambda_s:=\Diag (\lambda_1,...,\lambda_s)$. Denote $P=Q^\top \Diag (x)Q$, $P_s$ as the principle sub-matrix of $P$ indexed by $(1,...,s)$, $P_{n-s}$ as the principle sub-matrix of $P$ indexed by $(s+1,...,n)$ and $P_{s,n-s}$ as the sub-matrix of $P$ with rows indexed by $(1,...,s)$ and columns indexed by $(s+1,...,n)$ ($P_{n-s,s}$ similarly).
	
	Then the value $\lim_{\gamma\rightarrow +\infty}\linx (C,s; \gamma)$ exists and is the optimal value of the following convex program: 
	\begin{equation}\label{limcvx} 
	\begin{array}{cl}
	 \max &{\textstyle\frac{1}{2}} \!\left(\ldet \left(\Lambda_{s} P_{s} \Lambda_{s}\right)\!+\!\ldet \left(I_{n-s}-P_{n-s}\right)\! \right)\\ \text { s.t. } & \mathbf{e}^\top  x=s \\ & 0 \leq x \leq 1. 
	\end{array}
	\end{equation}
	Furthermore, $\linx (C,s; \gamma)$   is monotonically non-increasing in $\gamma$. 
\end{thm}

\proof{Proof.}

	By the conditions, 
	\begin{align*}
& 	\ldet (\gamma C \Diag (x)C+I-\Diag (x))-s\log \gamma\\
	& \quad =\ldet (\gamma \Lambda P \Lambda +I-P)-s\log\gamma.
	\end{align*} 
	Because $C,s$ are fixed, let $F_s(\gamma;x):=\gamma \Lambda_s P_s \Lambda_s +I_s-P_s$ be the principle sub-matrix of $\gamma \Lambda P \Lambda +I-P$ indexed by $(1,...,s)$.  We first prove that for any $\gamma>0$ and any $x$ feasible, $F_s(\gamma;x)$ is positive-definite so that we can use Schur complement formula to represent the determinant of $\gamma \Lambda P \Lambda +I-P$. 
	
	The construction of $P$ implies its eigenvalues are $\{x_1,x_2,...,x_n\}$ so all eigenvalues of $P$ lie in $[0,1]$. Because $P_s$ is a principle sub-matrix of $P$, by [Theorem 4.3.17, \cite{horn2012matrix}], all eigenvalues of $P_s$ lie in $[0,1]$. Decompose $P_s$ as $P_s=\hat Q \hat \Lambda \hat Q^\top $ where $\hat Q$ is orthogonal and $\hat \Lambda$ is the diagonal matrix of eigenvalues of $P_s$. In particular, all elements of $\diag(\hat \Lambda)$ are in $[0,1]$. Let $\hat C=\hat{Q}^\top  \Lambda_s\hat Q$, then 
	\begin{align*}
	F_s(\gamma;x)=\hat{Q}\left(\gamma \hat C\hat \Lambda \hat C+I_s-\hat \Lambda\right)\hat{Q}^\top .
	\end{align*}
	Because $\Lambda_s$ is positive-definite, so is $\hat C$. By [Lemma 20, \cite{chen_mixing}], $F_s(\gamma;x)$ is positive-definite for any $\hat\Lambda$ where $0\le \diag (\hat \Lambda)\le \mathbf{e}$. 
	
	
	We only need to consider $x$ in the feasible region such that $\gamma \Lambda P \Lambda +I-P$, (equivalently, $\gamma C \Diag (x)C+I-\Diag (x)$) is positive-definite. 
	So we assume that $\gamma \Lambda P \Lambda +I-P$ is positive-definite in the following. 
	Then the Schur complement of $\gamma \Lambda P \Lambda +I-P$ in $F_s(\gamma;x)$, which is $I_{n-s}-P_{n-s} -P_{n-s,s}F_s(\gamma;x)^{-1}P_{s,n-s}$, is also positive-definite. Furthermore, $P_{n-s,s}F_s(\gamma;x)^{-1}P_{s,n-s}$ is positive-semidefinite by the positive-definiteness of $F_s(\gamma;x)$ and we get that $I_{n-s}-P_{n-s}$ is positive-definite. On the other hand, because the feasible region of $x$ is compact, and the objective value of \eqref{limcvx} is upper bounded, the optimal value of \eqref{limcvx} is attainable by some $x$ such that the corresponding $\Lambda_s P_s \Lambda_s$ and $I_{n-s}-P_{n-s}$ are positive-definite. So we justify the definition of \eqref{limcvx}, and
	\begin{align*}
	& \ldet (\gamma \Lambda P \Lambda +I-P)-s\log\gamma\\
	&= \ldet (F_s(\gamma;x)) -s\log\gamma\\
&	+\!\ldet ( I_{n-s}\!-\!P_{n-s} \!-\!P_{n-s,s}F_s(\gamma;x)^{-1}\!P_{s,n-s} )\\
	&\textstyle = \ldet \left(\Lambda_s P_s \Lambda_s +\frac{1}{\gamma}(I_s-P_s)\right)\\
	&+\!\ldet (I_{n-s}\!-\!P_{n-s}\!-\!P_{n-s,s}F_s(\gamma;x)^{-1}\!P_{s,n-s})\!.
	\end{align*}
	Denote the optimal solution of \eqref{limcvx} as $x^*$ and $P^*=Q^\top \Diag(x^*) Q$. We claim that $\lim_{\gamma\rightarrow +\infty}\linx (C,s;\gamma)=$
	\begin{align*}
	{\textstyle\frac{1}{2}}\left( \ldet (\Lambda_s P^*_s \Lambda_s )+\ldet (I_{n-s}-P^*_{n-s})\right).
	\end{align*}
	
	We now prove this claim. In fact, for any $x$ feasible to \eqref{gammalinx} such that $\gamma \Lambda P \Lambda +I-P$ is positive-definite and that $F_s(\gamma;x)$ is positive-definite, we have
    \begin{align}\label{infub}
    &\textstyle  \ldet (\Lambda_s P_s \Lambda_s+\frac{1}{\gamma}(I_s-P_s))\\
    &\textstyle  +\ldet (I_{n-s}\!-\!P_{n-s}\!-\!P_{n-s,s}F_s(\gamma;x)^{-1}\!P_{s,n-s})\nonumber\\
 &  \textstyle   \le \ldet (\Lambda_s P_s \Lambda_s+\frac{1}{\gamma}I_s)+\ldet (I_{n-s}-P_{n-s}).\nonumber
    \end{align}
	
	We further assume that $\Lambda_s P_s \Lambda_s$ is positive-definite otherwise the right-hand-side of \eqref{infub} goes to minus infinity as $\gamma$ goes to infinity because $\ldet (I_{n-s}-P_{n-s})$ is clearly upper bounded by a uniform finite number for any $x$.
	
	Decompose $\Lambda_s P_s \Lambda_s$ as $Q'\Lambda'Q'^\top $ where $Q'$ is orthogonal and $\Lambda': =\Diag (\lambda_1',...,\lambda_s')$ where $\lambda_1'\ge \lambda_2'\ge\ldots\ge\lambda_s'>0$ is the diagonal matrix of eigenvalues of $\Lambda_s P_s \Lambda_s$. Then
	\begin{align*}
\textstyle	\ldet \left(\Lambda_s P_s \Lambda_s +\frac{1}{\gamma}I_s\right)=\log\left(\displaystyle\prod_{i=1}^s \textstyle \left(\lambda'_i+\frac{1}{\gamma}\right)\right).
	\end{align*}
	
Because every element of $\Lambda_s P_s \Lambda_s$ is bounded by a uniform number for any $x$, by Gershgorin circle theorem, $\lambda'_i,i\in\{1,...,s\}$ are bounded by a uniform number for all $x$. We pick a positive number $L_1>0$, when $\gamma\ge L_1$, there is a compact set $\mathcal{H}\subset \mathbb{R}^s$ (independent of $\gamma$) such that for all $x$ feasible to \eqref{gammalinx}, $(\lambda'_1+\frac{1}{\gamma},...,\lambda_s'+\frac{1}{\gamma})^\top$ as well as $(\lambda'_1,...,\lambda_s')^\top$ belongs to $\mathcal{H}$. Because the function $\prod_{i=1}^s y_i$ is continuous differentiable in $y$ on $\mathbb{R}^s$, it is Lipschitz continuous on $\mathcal{H}$, then, $\exists$ $L_2>0$ such that
	\begin{align*}
	\left|\prod_{i=1}^s \textstyle \left(\lambda'_i+\frac{1}{\gamma}\right)-\displaystyle \prod_{i=1}^s \lambda'_i\right|\leq \textstyle \frac{L_2\sqrt{s}}{\gamma}.
	\end{align*}
	Because $I_{n-s}-P_{n-s}$ is positive-definite and every element is bounded by a uniform number for any $x$, there exists $L_3>0$,
	\begin{align*}
	0<\det (I_{n-s}-P_{n-s})\le L_3.
	\end{align*}
	With the above arguments, when $\gamma\ge L_1$, we have 
	\begin{align*}
	&\textstyle \det\left(\Lambda_s P_s \Lambda_s +\frac{1}{\gamma}I_s\right)\det (I_{n-s}\!-\!P_{n-s})\\
	&= \left(\prod_{i=1}^s \textstyle\left(\lambda'_i+\frac{1}{\gamma}\right)\right)\det (I_{n-s}\!-\!P_{n-s})\\
	&\le  \left(\prod_{i=1}^s \lambda'_i+ \textstyle\frac{L_2\sqrt{s}}{\gamma}\right) \det (I_{n-s}\!-\!P_{n-s})\\
	& = \det (\Lambda_s P_s \Lambda_s)\!\det (I_{n-s}\!-\!P_{n-s})\\&\qquad +\textstyle\frac{L_2\sqrt{s}}{\gamma}\!\det (I_{n-s}\!-\!P_{n-s})\\
	& \le  \det (\Lambda_s P^*_s \Lambda_s)\!\det (I_{n-s}\!-\!P^*_{n-s})\\ &\qquad +\textstyle\frac{L_2\sqrt{s}}{\gamma}\!\det (\!I_{n-s}\!-\!P_{n-s})\\
	&\!\le\!   \det (\Lambda_s P^*_s \Lambda_s)\!\det (I_{n-s}\!-\!P^*_{n-s})\!+\!\textstyle\frac{L_2 L_3\sqrt{s}}{\gamma}.
	\end{align*}
	
	For any $x$ such that $\Lambda_s P_s \Lambda_s$ is singular, because the eigenvalues of $\Lambda_s P_s \Lambda_s$ are upper bounded uniformly for all $x$ feasible, clearly there is some $L_4>0$ such that when $\gamma\ge L_4$, any such $x$ cannot be an optimal solution for \eqref{gammalinx}.
	
	Because $\log(\cdot)$ is monotonically increasing, the above implies that when $\gamma\ge \max\{L_1,L_4\}$, we have
	\begin{align*}
	&\linx (C,s; \gamma)\\
	&=\max_{\begin{array}{c}
	 \scriptstyle     \mathbf{e}^\top x=s, \\
	 \scriptstyle     0\le x\le \mathbf{e}
	\end{array}} \textstyle {\textstyle\frac{1}{2}}\left(\ldet \left(\Lambda_s P_s \Lambda_s+\frac{1}{\gamma}(I_s-P_s)\right)\right.\\
	&\textstyle  \left. +
	\vphantom{\ldet \left(\Lambda_s P_s \Lambda+\frac{1}{\gamma}(I_s-P_s)\right)}
	\ldet (I_{n-s}-P_{n-s}-P_{n-s,s}F_s(\gamma;x)^{-1}P_{s,n-s})\right)\\
	 &\leq {\textstyle\frac{1}{2}}\log \left(\det (\Lambda_s P^*_s \Lambda_s)\det (I_{n-s}-P^*_{n-s})\right.\\
	 &\qquad \left.+\textstyle \frac{L_2 L_3\sqrt{s}}{\gamma}\right).
	\end{align*}
Taking limits on both sides, we have
\begin{align*}
&\lim_{\gamma\rightarrow +\infty}\linx (C,s; \gamma)\leq \\
&\lim_{\gamma\rightarrow +\infty} \!{\textstyle\frac{1}{2}}\!\log \!\left(\!\det (\Lambda_s P^*_s \Lambda_s)\!\det (I_{n-s}\!-\!P^*_{n-s})\right.\\
&\qquad  \left.+\textstyle \frac{L_2L_3\sqrt{s}}{\gamma}\right)\\
&= {\textstyle\frac{1}{2}}\left(\ldet (\Lambda_s P^*_s \Lambda_s )+\ldet (I_{n-s}-P^*_{n-s})\right).
\end{align*}
On the other hand, the optimal solution $x^*$ of \eqref{limcvx} is feasible to \eqref{gammalinx} and we have proved before that $\Lambda_s P^*_s \Lambda_s$ and $I_{n-s}-P^*_{n-s}$ are positive-definite, we have $\lim_{\gamma\rightarrow \infty} F_s(\gamma;x^*)^{-1}=O_{n}$ where $O_{n}$ is an all-zeros order-$n$ matrix and $\linx (C,s; \gamma)\ge f (C,s; \gamma; x^*)$. Furthermore,
\begin{align*}
& \lim_{\gamma\rightarrow +\infty}\linx (C,s; \gamma)
\ge  \lim_{\gamma\rightarrow +\infty}f (C,s; \gamma; x^*)\\
&= \lim_{\gamma\rightarrow +\infty}{\textstyle\frac{1}{2}}\left(\ldet \left(\Lambda_s P^*_s \Lambda_s \textstyle+\frac{1}{\gamma}(I_s-P^*_s)\right)\right.\\
& \left. +\ldet (I_{n-s}-P^*_{n-s}-P^*_{n-s,s}F_s(\gamma;x^*)^{-1}P^*_{s,n-s})\right)\\
&=  \textstyle\frac{1}{2}\left( \ldet (\Lambda_s P^*_s \Lambda_s )+\ldet (I_{n-s}-P^*_{n-s})\right).
\end{align*}

In all, we have
$\displaystyle \lim_{\gamma\rightarrow +\infty}\!\!\linx (C,s; \gamma)$
\[
=\!{\textstyle\frac{1}{2}}\!\!\left( \ldet (\Lambda_s P^*_s \Lambda_s )\!+\!\ldet (I_{n-s}\!-\!P^*_{n-s})\right)\!.
\]
Finally, because $\linx (C,s; e^{\psi})$ is convex in $\psi$ (see [Theorem 18, \cite{chen_mixing}])
and has a finite limit as $\psi\rightarrow +\infty$, we can conclude that $\linx (C,s; e^{\psi})$ is non-increasing in 
$\psi$, and hence $\linx (C,s; \gamma)$ is non-increasing in $\gamma$.
\Halmos
\endproof

At the outset, we assumed  $s\leq \rank(C)$.
Of course, the case where $s>\rank (C)$ is a bit strange because the optimal value of \ref{MESP} is always $-\infty$. But by the following theorem, the linx-bound problem can recognize these  cases.

\begin{thm}\label{slarrank}
If $s>\rank (C)$, then
$
\lim_{\gamma\rightarrow +\infty}
$
$
\linx  (C,s; \gamma)
$
$
=-\infty,
$
and there is no optimal $\gamma$. 
\end{thm}

\proof{Proof.}

Let $r=\rank (C)$. We use similar notations as in Theorem~\ref{seqrank}, but with the a little difference. Here we have $\Lambda_r: =\Diag\{\lambda_1,\lambda_2,\ldots,\lambda_r\}$ and  $P_{r}$, $P_{n-r}$, $P_{n-r,r}$, $P_{r,n-r}$ similarly because $r<s$. Then
\begin{align*}
&\ldet (\gamma C \Diag (x)C+I-\Diag (x))-s\log \gamma\\
&=\textstyle \ldet \left(\Lambda_r P_r \Lambda_r +\frac{1}{\gamma}(I_r-P_r)\right)-(s-r)\log\gamma\\
&\quad +\ldet (I_{n-r}-P_{n-r}-P_{n-r,r}F_r(\gamma;x)^{-1}P_{r,n-r}).
\end{align*} 

We consider the convex program
\begin{align} \label{rankcvxprogram}
&\textstyle \max~{\textstyle\frac{1}{2}}\left(\ldet \left(\Lambda_{r} P_{r} \Lambda_{r}  \right)\!+\!\ldet (I_{n-r}\!-\!P_{n-r})\!\right)\\
&\text{s.t.~} \mathbf{e}^\top x=s,~ 0\le x\le 1.\nonumber
\end{align}
Similar to that in Theorem~\ref{seqrank}, \eqref{rankcvxprogram} is well-defined. Denote the optimal solution of \eqref{rankcvxprogram} as $x^*$ and we have corresponding $P_r^*$ and $P_{n-r}^*$, by similar arguments as in Theorem~\ref{seqrank}, there exists $L_1, L_2, L_3, L_4 >0$ such that when $\gamma\ge \max\{L_1, L_4\}$,
\begin{align*}
&\lim_{\gamma\rightarrow +\infty} \linx  (C,s; \gamma)\\
& = \lim_{\gamma\rightarrow +\infty}\max_{\begin{array}{c}
	  \scriptstyle   \mathbf{e}^\top x=s, \\[-3pt]
	  \scriptstyle   0\le x\le \mathbf{e}
	\end{array}} {\textstyle\frac{1}{2}}\left(\ldet \left(\Lambda_r P_r \Lambda_r +\textstyle \frac{1}{\gamma}(I_r-P_r)\right)\right.\\[-3pt]
&\left. \qquad +\ldet (I_{n-r}-P_{n-r}-P_{n-r,r}F_r(\gamma;x)^{-1}P_{r,n-r})\right.\\[-3pt]
&\left.\vphantom{\ldet \left(\Lambda_r P_r \Lambda_r +\textstyle \frac{1}{\gamma}(I_r-P_r)\right)}\qquad  -(s-r)\log\gamma\right)\\[-3pt]
&\le  \lim_{\gamma\rightarrow +\infty}{\textstyle\frac{1}{2}} \left(\vphantom{\textstyle\frac{L_2 L_3 \sqrt{r}}{\gamma}}
\ldet (\Lambda_r P^*_r \Lambda_r )+\ldet (I_{n-r}-P^*_{n-r})\right.\\
& \left. \qquad +\textstyle\frac{L_2 L_3 \sqrt{r}}{\gamma}-(s-r)\log\gamma\right)
= -\infty. \Halmos
\end{align*}
\endproof

\section{Linear gap under optimal scaling}\label{sec:gammagap}
In Theorem~\ref{maskgap}, we constructed an infinite sequence $\{C_n\}_{n\in \mathcal{I}}$ where by choosing mask $I$, we decreased the linx bound by an amount that is at least linear in $n$ (specifically, $\approx .0312n$).
This is even the case when we choose optimal scaling parameters $\gamma$ (separately), with some sacrifice
in the constant. 

\begin{thm}\label{optgammalingap}
There is an infinite sequence of positive-semidefinite matrices $\{C_n\}_{n\in \mathcal{I}}$, such that 
\begin{align*}
\min_{\gamma>0}\linx \left(C_n,{\textstyle\frac{n}{2}};\gamma\right)-\min_{\bar \gamma>0}\linx \left(C_n,{\textstyle\frac{n}{2}}; I, \bar \gamma\right)\ge b n
\end{align*}
for some positive scalar $b\ge 0.024036$.
\end{thm}

\proof{Proof.}

We consider a crafted sequence of $C_n$. Assuming $n=4k$, and $C_n$ is block diagonal with $k$ blocks as $\left(\begin{array}{cc}
       1  & c_1 \\
       c_1  & 1
    \end{array}\right)$ and $k$ blocks as $\left(\begin{array}{cc}
       1  & c_2 \\
       c_2  & 1
    \end{array}\right)$ where $c_1\neq c_2$, $c_1^2\leq 1$, $c_2^2\leq 1$.

By Lemma~\ref{lb},
\begin{align*}
&\linx\left(C_n,{\textstyle\frac{n}{2}};  \gamma\right)\\
&\ge  \sum_{i=1}^{k}\left(\textstyle\frac{1}{2}\log\left(\frac{(1-c_1^2)^2}{4}\gamma+\frac{1+c_1^2}{2}+\frac{1}{4\gamma}\right)\right.\\
&\left.\qquad\qquad+\textstyle\frac{1}{2}\log \left(\frac{(1-c_2^2)^2}{4}\gamma+\frac{1+c_2^2}{2}+\frac{1}{4\gamma}\right)\right)\\
&=\textstyle k \left({\textstyle\frac{1}{2}}\log\left(\frac{(1-c_1^2)^2}{4}\gamma+\frac{1+c_1^2}{2}+\frac{1}{4\gamma}\right)\right.\\
&\qquad\qquad\left. + \textstyle\frac{1}{2}\log \left(\frac{(1-c_2^2)^2}{4}\gamma+\frac{1+c_2^2}{2}+\frac{1}{4\gamma}\right)\right).
\end{align*}
If $c_1^2, c_2^2<1$, the minimum of $\frac{(1-c_1^2)^2}{4}\gamma+\frac{1+c_1^2}{2}+\frac{1}{4\gamma}$ is $1$, achieved by the unique minimizer $\hat \gamma_1=1-c_1^2$,
and the minimum of $\frac{(1-c_2^2)^2}{4}\gamma+\frac{1+c_2^2}{2}+\frac{1}{4\gamma}$ is $1$, achieved by the unique minimizer $\hat \gamma_2=1-c_2^2$. If $c_1^2=1$, then no matter what value $\gamma$ is, $\frac{(1-c_2^2)^2}{4}\gamma+\frac{1+c_2^2}{2}+\frac{1}{4\gamma}$ is always greater than $1$. The case for $c_2^2=1$ is similar.

Thus we can choose $c_1^2\not= c_2^2$, then for all possible values of $c_1, c_2$,
\begin{align}\label{sec4lingap}
 b_{c_1,c_2} \!:= \min_{\gamma>0} {\textstyle} \sum_{i=1}^2  \textstyle  \log\!\left(\frac{(1-c_i^2)^2}{4}\gamma+\frac{1+c_i^2}{2}+\frac{1}{4\gamma}\right) \!> \!0. 
\end{align}
Then we have
\begin{align*} 
\min_{\gamma>0}\linx\left(C_n,{\textstyle\frac{n}{2}}; \gamma\right)\ge \textstyle\frac{1}{2} k b_{c_1,c_2}=\frac{b_{c_1,c_2}}{8} n.
\end{align*}
On the other hand, by Proposition~\ref{diagoptgamma}, we have $\min_{\bar \gamma}\linx\left(C_n,{\textstyle\frac{n}{2}};I,\bar \gamma\right)=\linx\left(C_n,{\textstyle\frac{n}{2}};I, 1\right)=0$.
So,
\begin{align*}
&\min_{\gamma}\linx\left(C_n,{\textstyle\frac{n}{2}}; \gamma\right)-\min_{\bar \gamma}\linx\left(C_n,{\textstyle\frac{n}{2}}; I,\bar \gamma\right)\\
&\qquad \ge 
\textstyle \frac{b_{c_1,c_2}}{8} n.
\end{align*}
Let $b=\frac{b_{c_1,c_2}}{8} $, we get what we want. In particular, if we set $c_1=0, c_2=1$, the optimal $\gamma$ for \eqref{sec4lingap} is $\hat \gamma=\frac{1+\sqrt{3}}{2}$ and $\frac{b_{0,1}}{8}=$
\begin{align*}
& \textstyle
\frac{1}{8}\!\left(\log\left(1\!+\!\frac{1}{2(1+\sqrt{3})}\right)\!+\!\log\left({\textstyle\frac{1}{2}}\!+\!\frac{1}{2(1+\sqrt{3})}\!+\!\frac{1+\sqrt{3}}{8}\right)\!\right)\\[-3pt]
&\qquad \ge 0.024036. \Halmos
\end{align*}
\endproof

\section{Final Remarks}\label{sec:final}
Our technical results establish the strong potential for
masking to improve on the (scaled) linx bound. 
So the next logical step is to work on optimizing the mask in this context. Similar work was carried out successfully for the spectral bound (see \cite{AnstreicherLee_Masked} and \cite{BurerLee}), where nonconvexity and nondifferentiability were the main difficulties to overcome. In the context of the linx bound, even at smooth points, it is not easy to get a handle on the necessary derivative information. There is also the potential to incorporate the ``mixing'' technique of \cite{chen_mixing} on top of mask optimization. 
We are currently working in this direction, and we plan to report on algorithmic results (with experimentation on benchmark data) for mask optimization in a future paper.

 \bibliographystyle{alpha}
 
 \bibliography{CFL_OperRes}

\end{document}